# A Mixed Integer Linear Programming Method for Dynamic Economic Dispatch with Valve Point Effect

Shanshan Pan [1,2], Jinbao Jian[2]† and Linfeng Yang[3]

**Abstract** – In this paper, a mixed integer linear programming (MILP) formulation is proposed to solve the dynamic economic dispatch with valve-point effect (DED-VPE). Based on piecewise linearization technique, the non-convex and non-smooth generation cost is reformulated into a linear lower approximation which is better than the quadratic one, yielding an MILP formulation for the DED-VPE. When the segment parameter is set appropriately, the MILP formulation can be solved by a mixed integer programming (MIP) solver directly and efficiently. Thus, a global optimal solution within a preset tolerance can be guaranteed for the MILP formulation. Simulation results show that the proposed MILP formulation can be solved to reliable solutions in reasonable time.

## 1. Introduction

Dynamic economic dispatch (DED) problem is an essential tool for real-time control of power system operation. To make the DED more accurate and practical, the valve-point effect (VPE) which makes the generation cost function non-convex and non-smooth should be considered. In order to address this problem, a large number of heuristic methods have been proposed, including genetic algorithm (GA), particle swarm optimization (PSO), differential evolution (DE), etc. A new list of heuristic methods applied to the dynamic economic dispatch with valve-point effect (DED-VPE) can be referred to [1]. However, heuristic methods do not provide an optimality gap so you have no clue how well of a solution you have obtained. Recently, an alternate method is to approximate the VPE cost via piecewise linearization technique [2], yielding an MIQP formulation for the DED-VPE. But when the MIQP formulation is directly solved by using a mixed integer programming (MIP) solver, the optimization will suffer convergence stagnancy and run out of memory. As a result, the multi-step method, the warm start technique and the range restriction scheme are required [2]. However, the range restriction scheme just restricts the solution space to a subspace where the global optimal solution would probably lie in. Consequently, the optimality of the solution for the MIQP can not be guaranteed.

In this paper, the DED-VPE is reformulated into an MILP which can obtain a better lower approximation for the generation cost in comparison with the MIQP. When the segment parameter is set appropriately, the MILP can be solved by a state-of-the-art MIP solver directly and efficiently. Thus, a global optimal solution within a preset tolerance can be guaranteed via an enumeration algorithm.

## 2. Mathematical Formulation of DED-VPE

The generation cost of each unit for conventional DED can be modeled by a convex quadratic polynomial (see Fig.1):

$$c_c(P_{i,t}) = \alpha_i + \beta_i P_{i,t} + \gamma_i P_{i,t}^2 \tag{1}$$

where $P_{i,t}$ is the power output of unit $i$ in period $t$; $\alpha_i$, $\beta_i$ and $\gamma_i$ are positive coefficients of unit $i$. When VPE is taken into account, a recurring rectified sinusoidal function (see Fig.1)

$$c_v(P_{i,t}) = |e_i \sin(f_i(P_{i,t} - P_i^{min}))| \tag{2}$$

is added to the conventional generation cost [3], which makes the generation cost function non-convex and non-smooth. Above-mentioned $P_i^{min}$ is the minimum power output of unit $i$; $e_i$ and $f_i$ are positive coefficients of unit $i$. Consequently, the generation cost for DED-VPE can be expressed as (see Fig.1):

$$c(P_{i,t}) = c_c(P_{i,t}) + c_v(P_{i,t}) \quad (3)$$

† Corresponding Author：jianjb@gxu.edu.cn
1  College of Electrical Engineering, Guangxi University, Nanning 530004, P. R. China.
2  Guangxi Colleges and Universities Key Laboratory of Complex System Optimization and Big Data Processing,Yulin Normal University, Yulin 537000, P. R. China.
3  College of Computer Electronics and Information, Guangxi University, Nanning 530004, P. R. China.

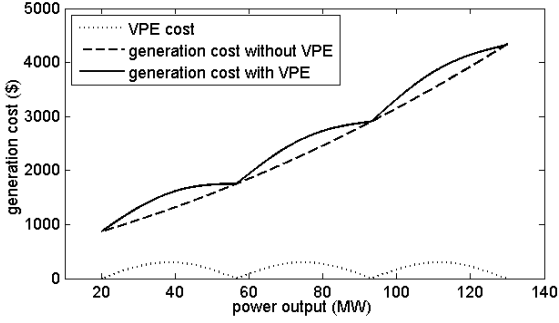

**Fig. 1.** The generation cost of unit for DED problem

The objective of DED-VPE is to minimize the total generation cost over a scheduled time horizon, which can be written as:

$$\min \sum_{t=1}^{T} \sum_{i=1}^{N} c(P_{i,t}) \quad (4)$$

where $N$ is the total number of units; $T$ is the total number of periods.

The minimized DED-VPE should be subjected to the constraints as follows.
1) Power balance equations

$$\sum_{i=1}^{N} P_{i,t} = D_t, \quad \forall\ t \quad (5)$$

where $D_t$ is the load demand in period $t$.
2) Power generation limits

$$P_i^{min} \leq P_{i,t} \leq P_i^{max}, \quad \forall\ i,t \quad (6)$$

where $P_i^{max}$ is the maximum power output of unit $i$.
3) Ramp rate limits

$$RD_i \leq P_{i,t} - P_{i,t-1} \leq RU_i, \quad \forall\ i,t \quad (7)$$

where $RD_i$ and $RU_i$ are the ramp-down and ramp-up rates of unit $i$, respectively.
4) Spinning reserve constraints

$$\begin{cases} SR_{i,t} \leq \min\{P_i^{max} - P_{i,t}, \tau RU_i\}, & \forall\ i,t \\ \sum_{i=1}^{N} SR_{i,t} \leq R_t, & \forall\ t \end{cases} \quad (8)$$

where $SR_{i,t}$ is the spinning reserve provided by unit $i$ in period $t$; $R_t$ is the system spinning reserve requirement in period $t$; $\tau$ is the time duration for units to deliver reserve [2].

### 3. Reformulation of DED-VPE

In [2], VPE cost (2) is considered for piecewise linearization and an MIQP is formed for DED-VPE. But when it is solved by an MIP solver directly, it suffers convergence stagnancy and will run out of memory, even for a 10-unit system. In other words, MIQP fails to address DED-VPE in a single step. Whereas, solution via MILP tends to be more efficient particularly because of the vastly superior warm start capabilities of the simplex method as compared with the interior-point one [4]. Therefore, different from [2], the whole generation cost (3) is considered for piecewise linearization in this paper.

#### 3.1 An MILP for DED-VPE

To obtain an MILP formulation of DED-VPE, $L_i + 1$ break points are chosen over a generation interval $[P_i^{min}, P_i^{max}]$, such that $P_i^{min} = a_{i,0} \leq a_{i,1} \leq \cdots \leq a_{i,L_i} = P_i^{max}$. Segment variables $P_{l,i,t}$ and binary variables $U_{l,i,t}$ ($l = 1, \cdots, L_i$) are introduced to make $P_{s,i,t} = P_{i,t}$ and $P_{l,i,t} = 0$ ($l \neq s$) when the $P_{i,t}$ lies in segment $s$ ($s \in \{1, \cdots, L_i\}$). Then generation cost $c(P_{i,t})$ can be approximately linearized as:

$$\hat{c}(P_{i,t}) = \sum_{l=1}^{L_i}(k_{l,i} P_{l,i,t} + b_{l,i} U_{l,i,t}) \quad (9)$$

with some additional constraints

$$\begin{cases} P_{i,t} = \sum_{l=1}^{L_i} P_{l,i,t} \\ a_{l-1} U_{l,i,t} \leq P_{l,i,t} \leq a_l U_{l,i,t} \\ \sum_{l=1}^{L_i} U_{l,i,t} = 1 \\ U_{l,i,t} \in \{0,1\}, \end{cases} \quad (10)$$

where $L_i$, $k_{l,i}$ and $b_{l,i}$ are calculated as follows

$$\begin{cases} L_i = ceil\left(M \dfrac{f_i\left(P_i^{max} - P_i^{min}\right)}{\pi}\right) \\ k_{l,i} = \dfrac{c(a_{l,i}) - c(a_{l-1,i})}{a_{l,i} - a_{l-1,i}} \\ b_{l,i} = c(a_{l-1,i}) - k_{l,i} a_{l-1,i}. \end{cases} \quad (11)$$

Foregoing $ceil(x)$ means round $x$ to the nearest integer greater than or equal to $x$ and $M$ is the number of equal segments on each $\sin(x)$ where $x$ belongs to $[0, \pi]$. It is well known that by refining the segments, i.e., by choosing $M$ large enough, an approximation of arbitrary accuracy can be achieved. However, the amount of extra continuous variables, binary variables and constraints for MILP formulation will increase significantly. For a tradeoff between modeling and computational efficiency, $M$ is usually not very large. In this paper, two cases: $M = 2$ and $M = 4$, will be adopted for numerical simulation. An example of piecewise linearization for a generation cost where $M = 2$ is shown in Fig. 2.

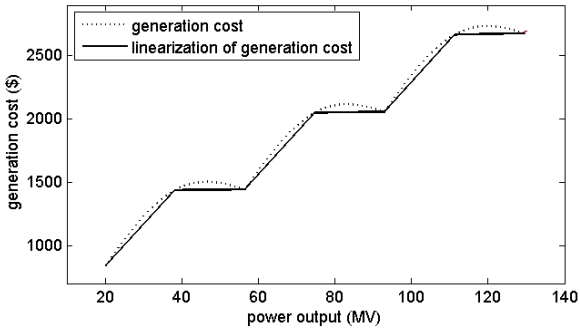

**Fig. 2.** Piecewise linearization for a generation cost ($M = 2$)

Consequently, the DED-VPE can be formulated as an MILP:

$$\min \sum_{t=1}^{T} \sum_{i=1}^{N} \hat{c}(P_{i,t}) \\ s.t. (5), (6), (7), (8), (10).$$
(12)

### 3.2 Comparison of Two Formulations

Generally, the quadratic approximation for a function is more accurate than the linear one. But this is not the case for the MIQP formulation [2] and our MILP formulation. Note that, although MIQP and MILP both can obtain a lower approximation for the generation cost function, but MILP can get a better lower approximation if the prespecified break points are the same. Since that, for any given two adjacent break points $a_{l-1,i}$, $a_{l,i}$, and for any $P_{i,t} \in (a_{l-1,i}, a_{l,i})$, by using the strict convexity of $c_c(P_{i,t})$, we have,

$$\hat{c}(P_{i,t}) = k_{l,i} P_{l,i,t} + b_{l,i}$$
$$= k_{l,i}(P_{l,i,t} - a_{l-1,i}) + c(a_{l-1,i})$$
$$= \frac{c(a_{l,i}) - c(a_{l-1,i})}{a_{l,i} - a_{l-1,i}}(P_{l,i,t} - a_{l-1,i}) + c(a_{l-1,i})$$
$$= \frac{c_c(a_{l,i}) - c_c(a_{l-1,i})}{a_{l,i} - a_{l-1,i}}(P_{l,i,t} - a_{l-1,i}) + c_c(a_{l-1,i})$$
$$+ \frac{c_v(a_{l,i}) - c_v(a_{l-1,i})}{a_{l,i} - a_{l-1,i}}(P_{l,i,t} - a_{l-1,i}) + c_v(a_{l-1,i})$$
$$= \frac{c_c(a_{l,i}) - c_c(a_{l-1,i})}{a_{l,i} - a_{l-1,i}}(P_{l,i,t} - a_{l-1,i}) + c_c(a_{l-1,i}) + \hat{c}_v(P_{i,t})$$
$$= \alpha c_c(a_{l,i}) + \beta c_c(a_{l-1,i}) + \hat{c}_v(P_{i,t})$$
$$\geq c_c(P_{i,t}) + \hat{c}_v(P_{i,t})$$
(13)

where

$$\begin{cases} \alpha = \frac{P_{l,i,t} - a_{l-1,i}}{a_{l,i} - a_{l-1,i}} \\ \beta = 1 - \frac{P_{l,i,t} - a_{l-1,i}}{a_{l,i} - a_{l-1,i}} = \frac{a_{l,i} - P_{l,i,t}}{a_{l,i} - a_{l-1,i}} \end{cases}$$
(14)

and $\hat{c}_v(P_{i,t})$ is the linear formulation of $c_v(P_{i,t})$.

### 3.3 The Estimation of Optimality

As we know, our MILP formulation can obtain a lower approximation for the generation cost function. Thus, its best lower bound within a relative mipgap tolerance *RGap* is also a lower bound for the original DED-VPE. Then the optimality gap of the original problem, denoted by *OGap*, can be defined as:

$$OGap = (Z - Lb)/Lb$$
(15)

where $Z$ is the obtaining optimal value of the original problem, $Lb$ is the corresponding lower bound provided by an MIP solver.

Thereby, the quality of the solution we gained can be measured via *OGap*.

## 4. Simulation Results

In this section, a set of different sizes test systems with units ranging from 10 to 500 over a scheduled time horizon of 24 h and two cases: $M = 2$ and $M = 4$, are adopted for testing the effectiveness of the proposed MILP formulation. The 10-unit system is taken from [5]. The 30-, 100- and 500-unit systems are obtained by duplicating the 10-unit system three, ten and fifty times. For fair comparison, the 1-h spinning reserve requirement is 5% of the load demand and the 10-min spinning reserve requirement is $(2/6) \times 5\%$ of the load demand. Meanwhile, the scaled CPU time [6]:

$$Scaled\ CPU\ time = \frac{Given\ CPU\ speed}{Base\ CPU\ speed} Given\ CPU\ time$$
(16)

is used in this paper and the base CPU speed is 2.4 GHz.

The model is coded with Matlab and optimized using Cplex 12.6.2. The machine for all runs is an Intel Core 2.5 GHz Dell-notebook with 8 GB of RAM.

The simulation results obtained by MILP are listed in Table 1, where the costs are calculated by objective function (4) to eliminate the error caused by piecewise linearization.

According to Table 1, our MILP formulation out performs the MIQP formulation in terms of the total generation cost. When $M = 2$ and $RGap$ is set to 0.25%, the 10-, 30- and 100-unit systems all can be solved to lower costs in a faster speed. For the 500-unit system, more time is consumed than MIQP. Actually, it is reasonable that more time is expended, since our MILP formulation is solved to a much lower cost by Cplex directly. When a smaller $RGap$ 0.20% is set, more superior results can be found and of course, more time is required.

As we know, a large $M$ which makes the MILP formulation more accurate will result in a better solution. But at the same time, the computational efficiency may greatly reduce. For the 30- and 500-unit systems, when $RGap = 0.25\%$, the case $M = 4$ can solve to more accurate solutions than the case $M = 2$. Nevertheless, with a smaller $RGap$ 0.20% setting, the case $M = 2$ can get lower costs in shorter time in comparison with the case $M = 4$. For the 100-unit system, a larger cost is obtained in the case $M = 4$ when $RGap = 0.25\%$. This abnormal performance mainly because when $M = 2$, the same cost 10154980 $ is calculated with $RGap$ 0.22%~0.30%. Indeed, when the $RGap$ is set to 0.20% for the case $M = 4$, a cost 10150983 $ which is smaller than the cost 10151410 $ for the case $M = 2$ can be obtained by using a much longer time 13.74 min. But for a small 10-unit system, the case $M = 4$ seems to be more efficiency.

Since the 30-, 100- and 500-unit systems are obtained by duplicating the 10-unit system. Thus, the average cost, denoted as A-cost, for per ten units will be no more than the 10-unit system cost. When $M = 2$, all the A-costs obtained by MILP are lower than the 10-unit system cost, which conforms to the accuracy of the results. Whereas, in MIQP, the A-costs for the 100-unit and 500-unit systems are much larger than the corresponding 10-unit system cost, which implies that the optimality of the MIQP can not be guaranteed. Actually, when $M = 1$, the cost for the 500-unit system obtained by MILP within 0.42 min is 50834169 $, which is much lower than the cost 51354130 $ obtained by MIQP.

Simultaneously, we can see that all the $OGaps$ for our solutions are no more than 0.5%, which indicates that our MILP formulation can be solved to more reliable and efficient solutions instead of probably global optimal solutions.

The output for each unit obtained by the MILP for the 10-unit system is given in Table 2 for verification. The outputs of unit 7 and 10 are not listed in the table because they are always 129.59 MW and 55.00 MW, respectively.

## 5. Conclusion

**Table 1.** Results for the 10-, 30-, 100- and 500-unit systems

| System | Method | $M$ | Cost ($) | Time (min) | A-Cost ($) | $RGap$ | $OGap$ |
|---|---|---|---|---|---|---|---|
| 10-unit | MIQP [2] | / | 1016601 | 1.88 | / | / | / |
| | MILP | 2 | 1016533 | 0.50 | / | 0.25% | 0.40% |
| | | | 1016429 | 3.80 | / | 0.20% | 0.36% |
| | | 4 | 1016329 | 3.80 | / | 0.30% | 0.32% |
| 30-unit | MIQP [2] | / | 3049359 | 3.86 | 1016453 | / | / |
| | MILP | 2 | 3046454 | 0.20 | 1015485 | 0.25% | 0.38% |
| | | | 3045922 | 0.27 | 1015307 | 0.20% | 0.35% |
| | | 4 | 3046135 | 1.68 | 1015378 | 0.25% | 0.23% |
| 100-unit | MIQP [2] | / | 10170508 | 3.64 | 1017051 | / | / |
| | MILP | 2 | 10154980 | 1.38 | 1015498 | 0.25% | 0.41% |
| | | | 10151410 | 3.15 | 1015141 | 0.20% | 0.30% |
| | | 4 | 10155601 | 10.93 | 1015560 | 0.25% | 0.25% |
| 500-unit | MIQP [2] | / | 51354130 | 12.83 | 1027083 | / | / |
| | MILP | 2 | 50773342 | 66.70 | 1015467 | 0.25% | 0.34% |
| | | | 50749042 | 72.93 | 1014981 | 0.20% | 0.30% |
| | | 4 | 50759312 | 243.46 | 1015186 | 0.25% | 0.22% |

**Table 2.** Outputs (MW) of units for the 10-unit system ($M = 2, RGap = 25\%$)

| $t$ | Unit 1 | Unit 2 | Unit 3 | Unit 4 | Unit 5 | Unit 6 | Unit 8 | Unit 9 |
|---|---|---|---|---|---|---|---|---|
| 1 | 150.00 | 222.27 | 156.69 | 60.00 | 73.00 | 122.45 | 47.00 | 20.00 |
| 2 | 150.00 | 229.53 | 223.43 | 60.00 | 73.00 | 122.45 | 47.00 | 20.00 |
| 3 | 150.00 | 309.53 | 291.43 | 60.00 | 73.00 | 122.45 | 47.00 | 20.00 |
| 4 | 150.00 | 389.53 | 297.40 | 60.00 | 122.87 | 134.61 | 47.00 | 20.00 |
| 5 | 226.62 | 396.80 | 297.40 | 60.00 | 122.87 | 124.72 | 47.00 | 20.00 |
| 6 | 303.25 | 396.80 | 297.40 | 60.00 | 172.73 | 146.23 | 47.00 | 20.00 |
| 7 | 379.87 | 396.80 | 297.40 | 70.42 | 172.73 | 133.19 | 47.00 | 20.00 |
| 8 | 379.87 | 396.80 | 282.27 | 120.42 | 222.60 | 122.45 | 47.00 | 20.00 |
| 9 | 456.50 | 396.80 | 297.40 | 170.42 | 222.60 | 122.45 | 53.25 | 20.00 |
| 10 | 456.50 | 396.80 | 297.85 | 220.42 | 222.60 | 160.00 | 83.25 | 50.00 |
| 11 | 456.50 | 396.80 | 340.00 | 248.14 | 222.60 | 160.00 | 85.31 | 52.06 |
| 12 | 456.50 | 460.00 | 327.70 | 241.25 | 222.60 | 160.00 | 115.31 | 52.06 |
| 13 | 456.50 | 396.80 | 297.40 | 218.80 | 222.60 | 160.00 | 85.31 | 50.00 |
| 14 | 456.50 | 396.80 | 296.95 | 168.80 | 222.60 | 122.45 | 55.31 | 20.00 |
| 15 | 379.87 | 396.80 | 297.40 | 118.80 | 172.73 | 158.80 | 47.00 | 20.00 |
| 16 | 303.25 | 396.80 | 288.24 | 68.80 | 122.87 | 122.45 | 47.00 | 20.00 |
| 17 | 226.62 | 396.80 | 297.40 | 60.00 | 122.87 | 124.72 | 47.00 | 20.00 |
| 18 | 303.25 | 396.80 | 297.40 | 70.42 | 172.73 | 131.07 | 51.74 | 20.00 |
| 19 | 379.87 | 396.80 | 297.40 | 120.42 | 172.73 | 122.45 | 81.74 | 20.00 |
| 20 | 456.50 | 402.45 | 340.00 | 170.42 | 222.73 | 160.00 | 85.31 | 50.00 |
| 21 | 456.50 | 389.53 | 322.60 | 120.42 | 222.60 | 122.45 | 85.31 | 20.00 |
| 22 | 379.87 | 309.53 | 283.09 | 70.42 | 172.73 | 122.45 | 85.31 | 20.00 |
| 23 | 303.25 | 229.53 | 204.00 | 60.00 | 122.87 | 122.45 | 85.31 | 20.00 |
| 24 | 226.62 | 222.27 | 189.76 | 60.00 | 73.00 | 122.45 | 85.31 | 20.00 |

(*The outputs of unit 7 and 10 are always 129.59 MW and 55.00 MW.)

In this paper, a good MILP formulation is proposed for DED with VPE. When the segments are chosen appropriately, our MILP formulation can be solved to global optimality within a preset tolerance in a high efficiency. Simulation results demonstrate that the proposed MILP formulation is a promising tool for solving a practical DED problem.

## Acknowledgements

This work was supported by the Natural Science Foundation of China (No.51407037); the Natural Science Foundation of Guangxi (No.2014GXNSFFA118001); the Open Project Program of Guangxi Colleges and Universities Key Laboratory of Complex System Optimization and Big Data Processing (No.2016CSOBDP0205).